\documentclass[10pt]{amsart}
\usepackage{amsmath}
\usepackage{amssymb}

\def\underset#1#2{{\mathrel{\mathop {{}_{} {#2}}\limits_{{#1}_{}}}}}
\def\upplim_#1{\underset{#1}{\overline\lim}\;}
\def\lowlim_#1{\underset{#1}{\underline\lim}\;}
\newtheorem{corollary}[equation]{Corollary}

\newtheorem{claim}[equation]{Claim}

\newtheorem{lemma}[equation]{Lemma}

\newtheorem{theorem}[equation]{Theorem}

\newcommand{\C}{{\mathbb{C}}}

\newcommand{\N}{\mathbf{N}}

\renewcommand{\P}{{\mathbb{P}}}
\newcommand{\bbf}{{\bf{f}}}
\newcommand{\bbF}{{\bf{F}}}

\newcommand{\supp}{\mathrm{Supp}\,}

\newcommand{\Z}{\mathbb{Z}}

\numberwithin{equation}{section}

\title[Defect relation for holomorphic maps from complex discs]{Defect relation for holomorphic maps from complex discs into projective varieties and hypersurfaces} 

\author{Si Duc Quang}
\address{Department of Mathematics, Hanoi National University of Education, 136-Xuan Thuy, Cau Giay, Hanoi, Vietnam}
\email{quangsd@hnue.edu.vn}
\begin{document}

\maketitle 

\begin{abstract}
In this paper, we establish a second main theorem for holomorphic maps with finite growth index on complex discs intersecting arbitrary families of hypersurfaces (fixed and moving) in projective varieties, which gives an above bound of the sum of truncated defects. Our result also is generalizes and improves many previous second main theorems for holomorphic maps from $\C$ intersecting hypersurfaces (moving and fixed) in projective varieties.
\end{abstract}

\def\thefootnote{\empty}

\section{Introduction}

Let $\Delta(R)=\{z\in\C: |z|<R\}$ be a complex disc and $r_0$ be a fixed positive number so that $0<r_0<R$. Let $\nu$ be a divisor on $\Delta (R)$, which is regarded as a function on $\Delta (R)$ with values in $\Z$ such that $\supp\nu:=\{z;\nu(z)\ne 0\}$ is a discrete subset of $\Delta (R)$. Let $k$ be a positive integer or $+\infty$. The truncated counting function of $\nu$ is defined by:
$$ n^{[k]}(t)=\sum_{|z|\le t}\min\{k,\nu(z)\} \ (0\le t\le R),$$
$$ \text{ and }\ N^{[k]}(r,\nu)=\int_{r_0}^{r}\dfrac{n^{[k]}(t)-n^{[k]}(0)}{t}dt.$$
We will omit the character $^{[k]}$ if $k=+\infty$.

Let $\varphi:\Delta (R)\rightarrow \C\cup\{\infty\}$ be a non-constant meromorphic function. We denote by $\nu^0_\varphi$ (resp. $\nu^\infty_\varphi$) the divisor of zeros (resp. divisor of poles) of $\varphi$ and set $\nu_\varphi=\nu^0_\varphi-\nu^\infty_\varphi$. As usual, we will write $N^{[k]}_{\varphi}(r)$ and $N^{[k]}_{1/\varphi}(r)$ for $N^{[k]}(r,\nu^0_{\varphi})$ and $N^{[k]}(r,\nu^\infty_{\varphi})$ respectively.

Let $f:\Delta (R)\rightarrow\P^N(\C)$ be a holomorphic map and $\Omega$ be the Fubini-Study form on $\P^N(\C)$. The characteristic function of $f$ is defined by
$$ T_f(r):=\int_{0}^r\frac{dt}{t}\int_{|z|<t}f^*\Omega.$$
In \cite{RS}, M. Ru and N. Sibony defined the growth index of $f$ by
$$ c_{f}=\mathrm{inf}\left\{c>0\ \biggl |\int_{0}^R\mathrm{exp}(cT_{f}(r))dr=+\infty\right\}.$$
For the convenient, we will set $c_f=+\infty$ if 
$$\left\{c>0\ \bigl |\int_{0}^R\mathrm{exp}(cT_{f}(r))dr=+\infty\right\}=\emptyset.$$ 

A meromorphic function $a$ on $\Delta (R)$ (which is regarded as a holomorphic map into $\P^1(\C)$) is said to be small with respect to $f$ if $\| T_a(r)=o(T_f(r))$ as $r\rightarrow R$. Here (and throughout this paper), the notation ``$\| P$'' means the assertion $P$ holds for all $r\in (0,R)$ outside a subset $S$ of $(0,R)$ with $\int_{S}{\rm exp}((c_f+\epsilon)T_f(r))dr<+\infty$ for some $\epsilon>0$. 

Denote by $\mathcal H$ the ring of all holomorphic functions on $\Delta (R)$. Let $Q$ be a homogeneous polynomial in $\mathcal H[x_0,\dots,x_n]$ of
degree $d \geq 1$ given by
$$ Q(z)=\sum_{I\in\mathcal T_d}a_I(z)\omega^I, $$
where $\mathcal T_d=\{(i_0,\ldots,i_N)\in\N_0^{N+1}\ ;\ i_0+\cdots +i_N=d\}$, $\omega^I=\omega_0^{i_0}\cdots\omega_n^{i_N}$ for $I=(i_0,\ldots,i_N)$ and all $a_I\in\mathcal H$ has no common zero. The homogeneous polynomial  $Q$ is called a moving  hypersurface of $\P^N (\mathbf{C} )$. Throughout this paper, by changing the homogeneous coordinates of $\P^N(\C)$ if necessary, we may assume that $a_{I_0}\not\equiv 0$ each such given moving hypersurface $Q$, where $I_0=(d,0,\ldots,0)$. We put $\tilde Q=\sum_{I\in\mathcal T_d}\frac{a_{I}}{a_{I_0}}\omega^I$. 

The moving hypersurface $Q$ is said to be slow with respect to $f$ all $\frac{a_I}{a_{I_0}}\ (I\in\mathcal T_d)$ are small with respect to $f$.  Let $\bbf =(f_0,\ldots,f_n)$ be a reduced representation of $f$. We define 
$$ Q(\bbf)(z)=\sum_{I=(i_0,\ldots,i_n)\in\mathcal T_d}a_I(z)f_0^{i_0}(z)\cdots f_n^{i_n}(z).$$
Then the truncated divisor $\nu^{[k]}_{Q(\bbf)}$ does not depend on the choice of the reduced representation $\bbf$ and hence is written by $\nu_{Q(f)}$. Its truncated counting function is denoted simply by $N^{[k]}_{Q(f)}(r)$. The proximity function of $f$ with respect to $Q$ is define by
$$ m_f(r,Q)=\int_{0}^{2\pi}\log\frac{\|f\|^q\cdot\|Q\|}{|Q(f)|}(re^{i\theta})\frac{d\theta}{2\pi}.$$
If $Q$ is slow with respect to $f$, then the first main theorem states that
$$\|\ dT_f(r)=m_f(r,Q)+N_{Q(f)}(r)+o(T_f(r)).$$
The truncated defect of $f$ with respect to $Q$ is defined by
$$ \delta^{[k]}_{f,Q}= 1-\underset{r\longrightarrow R}{\rm limsup}\frac{N^{[k]}_{Q(f)}(r)}{d T_f(r)}.$$
We omit the character $^{[k]}$ if $k=+\infty$. If all coefficients of $Q$ are constant then we call $Q$ a (fixed) hypersurface of $\P^N(\C)$ and set $Q^*=\{(\omega_0:\cdots:\omega_n)\ |\ \sum_{I\in\mathcal T_d}a_I\omega^I=0\}$. 

Let $V$ be a smooth projective subvariety of $\P^N(\C)$ of dimension $n$. Let $\mathcal Q=\{Q_1,\ldots,Q_q\}$ be a family of  moving hypersurfaces in $\P^N(\mathbf{C})$, where $Q_i=\sum_{I\in\mathcal T_{d_i}}a_{iI}x^I$. Denote by $\mathcal K_{\mathcal Q}$ the smallest field which contains $\C$ and all functions $\frac{a_{iI}}{a_{iI_0}}\ (I\in\mathcal T_{d_i})$. As usual, the family $\{Q_1,\ldots,Q_q\}$ is said to be in weakly $\ell$-subgeneral position if $\bigcap_{s=1}^{\ell+1}Q_{j_s}(z)^*\cap V=\emptyset$ for every $1\le j_1<\cdots<j_{\ell+1}\le q$ and for generic points $z\in\Delta(R)$ (i.e., for all $z\in\Delta(R)$ outside a discrete subset). Here, we note that $\dim\emptyset=-\infty$. In \cite{Q22c}, we define the distributive constant of the family $\mathcal Q$ with respect to $V$ by
$$ \Delta_{V}:=\underset{\Gamma\subset\{1,\ldots,q\}}\max\dfrac{\sharp\Gamma}{n-\dim\left (\bigcap_{j\in\Gamma} Q_j(z)^*\right)}$$
for generic points $z\in\Delta (R)$. 
From \cite[Remark 3.7]{Q22c}, we know that if $\mathcal Q$ is in weakly $\ell$-subgeneral position with respect to $V$ then $\Delta_{V}\le\ell-n+1$.

For the case of holomorphic curve from $\C$ into $V$ and family of fixed hypersurfaces in general position (i.e., in $n$-subgeneral position), M. Ru \cite{Ru09} proved the following.

\vskip0.2cm
\noindent
{\bf Theorem A.} (see \cite{Ru09}) {\it Let $f$ be an algebraically nondegenerate holomorphic map of $\C$ into a smooth subvariety $V\subset\P^N(\mathbf{C})$ of dimension $n$. Let $\{Q_i\}_{i=1}^q$ be a family of $q$ hypersurfaces in general position with respective to $V$. Then for any $\epsilon >0$, we have
$$\|\ (q-n-1-\epsilon)T_f(r)\le \sum_{i=1}^{q}\dfrac{1}{\deg Q_i}N_{Q_i(f)}(r).$$}

\vskip0.1cm
\noindent
From the above theorem, the number $n+1$ is an above bound of the sum of defects (without truncated multiplicity) for hypersufaces in this case. Later on, many mathematicians generalized Theorem A to the case of slowly moving hypersurfaces in general position with respect to $V$. Recently, in \cite{Q22a} the author considered the case of holomorphic maps from $\Delta(R)$ into $\P^N(\C)$  and proved the following result.

\noindent
{\bf Theorem B.} (reformulation of \cite[Theorem 1.3]{Q22a}) {\it Let $V\subset\P^N(\C)$ be a smooth complex projective variety of dimension $n\ge 1$. Let $\{Q_1,\ldots,Q_q\}$ be a family of hypersurfaces in $\P^N(\C)$ with the distributive constant $\Delta$ with respect to $V$, $\deg Q_i=d_i\ (1\le i\le q)$, and let $d$ be the least common multiple of $d_1,\ldots,d_q$. Let $f:\Delta(R)\to V$ be an algebraically non-degenerate holomorphic curve with $c_f<+\infty$. Then, for every $\epsilon >0$, 
\begin{align*}
\bigl\|\ &(q-\Delta(n+1)-\epsilon)T_f(r)\\
&\le\sum_{i=1}^q\frac{1}{d_i}N^{[L]}_{Q_i(f)}(r)+\dfrac{Lc_fT_f(r)}{2d^{n+1}(2n+1)(n+1)(q!)\deg V}.
\end{align*}
where $L=\bigl[d^{n^2+n}\deg (V)^{n+1}e^n\Delta^n(2n+4)^n(n+1)^n(q!)^n\epsilon^{-n}\bigl].$}

Here, by $[x]$ we denote the largest integer not exceeding the real number $x$. However, the coefficient of $T_f(r)$ in the right hand side of the above inequality has a factor $c_f(q!\epsilon^{-1})^{n-1}$. Then if $c_f$ or $q$ is large enough, this coefficient always exceeds $q$ and the theorem is meaningless. Hence, it may not imply the defect relation. Our purpose in this paper is to improve Theorem B by reducing the truncation level $L$ and that coefficients so that they do not depend on $q$. In order to do so, we will apply the new below bound of Chow weight in \cite{Q22b} (see Lemma \ref{2.2}) and also give some new technique to control the error term occuring when the theorem on the estimate of Hilbert weights (Theorem \ref{2.1}) is applied. We also consider the case of moving hypersurfaces. Our main result is stated as follows.

\begin{theorem}\label{1.1} 
Let $f$ be a nonconstant holomorphic map of $\Delta (R)$ into an $n$-dimension smooth projective subvariety $V\subset\P^N(\mathbf{C})$ with finite growth index $c_f$. Let $\{Q_i\}_{i=1}^q$ be a family slow (with respect to $f$) moving hypersurfaces with the distributive constant $\Delta_{V}$ with respect to $V$. Assume that $f$ is algebraically nondegenerate over $\mathcal K_{\mathcal Q}$. 

a) For any $\epsilon'>0$ and $(n+1)\Delta_V>\epsilon>0$, we have
\begin{align*}
\bigl\|\ &(q-\Delta_{V}(n+1)-\epsilon)T_f(r)\\
&\le\sum_{j=1}^{q}\frac{1}{d}N^{[L-1]}_{Q_j(f)}(r)+\frac{(\Delta_{V}(n+1)+\epsilon)(c_f+\epsilon')(L-1)}{2du}T_f (r),
\end{align*}
where 
$$L=d^n\deg V(u+1)^n\left[\biggl(1+\frac{\epsilon}{2(n+1)\Delta_{V}}\biggl)^{\left[\frac{d^n\deg V(u+1)^{n+q}}{\log^2(1+\frac{\epsilon}{2(n+1)\Delta_{V}})}\right]+1}\right]$$
with $u=\lceil 2\Delta_{V}(2n+1)(n+1)d^n\deg V(\Delta_{V}(n+1)+\epsilon)\epsilon^{-1}\rceil$.

b) Assume further that all $Q_i\ (1\le i\le q)$ are assumed to be fixed hypersurfaces. For any $\epsilon'>0$ and $\epsilon>0$ we have 
\begin{align*}
\bigl\|\ &(q-\Delta_{V}(n+1)-\epsilon)T_f(r)\\
&\le\sum_{j=1}^{q}\frac{1}{d}N^{[L'-1]}_{Q_j(f)}(r)+\frac{(\Delta_{V}(n+1)+\epsilon)(c_f+\epsilon')(L'-1)}{2du'}T_f (r),
\end{align*}
where $L'=\bigl[d^{n^2+n}(\deg V)^{n+1}e^n(2n+5)^n(\Delta^2_V(n+1)\epsilon^{-1}+\Delta_V)^n\bigl]$ with $u'=\lceil \Delta_{V}(2n+1)(n+1)d^n\deg V(\Delta_{V}(n+1)+\epsilon)\epsilon^{-1}\rceil$.
\end{theorem}
Here, $\lceil x\rceil$ stands for the smallest integer bigger than or equal to the real number $x$.
By this theorem, we get the following truncated defect relation for fixed hypersurfaces.
\begin{corollary}\label{1.3}
With the assumption and notation as in Theorem \ref{1.1} and suppose that all $Q_i$ are fixed hypersurfaces. Then for any $\epsilon>0$, we have
 $$\sum_{i=1}^q\delta^{[L-1]}_{f,Q_i}\le\Delta_{V}(n+1)+\epsilon+\frac{(\Delta_{V}(n+1)+\epsilon)c_f(L-1)}{2du},$$
where $L=[d^{n^2+n}(\deg V)^{n+1}e^n(2n+5)^n(\Delta^2_V(n+1)\epsilon^{-1}+\Delta_V)^n]$ with $u=\lceil 2\Delta_{V}(2n+1)(n+1)d^n\deg V(\Delta_{V}(n+1)+\epsilon)\epsilon^{-1}\rceil$.
\end{corollary}

\section{Some auxiliary results}

Let $X\subset\P^n(\C)$ be a projective variety of dimension $k$ and degree $\delta$. For $\textbf{a} = (a_0,\ldots,a_n)\in\mathbb Z^{n+1}_{\ge 0}$ we write ${\bf x}^{\bf a}$ for the monomial $x^{a_0}_0\cdots x^{a_n}_n$. Let $I_X$ be the prime ideal in $\C[x_0,\ldots,x_n]$ defining $X$. Let $\C[x_0,\ldots,x_n]_u$ be the vector space of homogeneous polynomials in $\C[x_0,\ldots,x_n]$ of degree $u$ (including $0$). For $u= 1, 2,\ldots,$ put $(I_X)_u:=\C[x_0,\ldots,x_n]_u\cap I_X$ and define the Hilbert function $H_X$ of $X$ by
\begin{align*}
H_X(u):=\dim\C[x_0,\ldots,x_n]_u/(I_X)_u.
\end{align*}
Let ${\bf c}=(c_0,\ldots,c_n)$ be a tuple in $\mathbb R^{n+1}_{\ge 0}$ and let $e_X({\bf c})$ be the Chow weight of $X$ with respect to ${\bf c}$. The $u$-th Hilbert weight $S_X(u,{\bf c})$ of $X$ with respect to ${\bf c}$ is defined by
\begin{align*}
S_X(u,{\bf c}):=\max\sum_{i=1}^{H_X(u)}{\bf a}_i\cdot{\bf c},
\end{align*}
where the maximum is taken over all sets of monomials ${\bf x}^{{\bf a}_1},\ldots,{\bf x}^{{\bf a}_{H_X(u)}}$ whose residue classes modulo $I_X$ form a basis of $\C[x_0,\ldots,x_n]_u/(I_X)_u.$

The following theorem is due to J. Evertse and R. Ferretti \cite{EF01}.
\begin{theorem}[{see \cite[Theorem 4.1]{EF01}}]\label{2.1}
Let $X\subset\P^n(\C)$ be an algebraic variety of dimension $k$ and degree $\delta$. Let $u>\delta$ be an integer and let ${\bf c}=(c_0,\ldots,c_n)\in\mathbb R^{n+1}_{\geqslant 0}$.
Then
$$ \dfrac{1}{uH_X(u)}S_X(u,{\bf c})\ge\dfrac{1}{(k+1)\delta}e_X({\bf c})-\dfrac{(2k+1)\delta}{u}\cdot\left (\max_{i=0,\ldots,n}c_i\right).$$
\end{theorem}

The following lemma is due to the author \cite{Q22b}.
\begin{lemma}[{see \cite[Lemma 3.2]{Q22b}}]\label{2.2}
Let $Y$ be a projective subvariety of $\P^R(\C)$ of dimension $k\ge 1$ and degree $\delta_Y$. Let $\ell\ (\ell\ge k+1)$ be an integer and let ${\bf c}=(c_0,\ldots,c_R)$ be a tuple of non-negative reals. Let $\mathcal H=\{H_0,\ldots,H_R\}$ be a set of hyperplanes in $\P^R(\C)$ defined by $H_{i}=\{y_{i}=0\}\ (0\le i\le R)$. Let $\{i_1,\ldots, i_\ell\}$ be a subset of $\{0,\ldots,R\}$ such that:
\begin{itemize}
\item[(1)] $c_{i_\ell}=\min\{c_{i_1},\ldots,c_{i_\ell}\}$,
\item[(2)] $Y\cap\bigcap_{j=1}^{\ell-1}H_{i_j}\ne \emptyset$, 
\item[(3)] and $Y\not\subset H_{i_j}$ for all $j=1,\ldots,\ell$.
\end{itemize}
Let $\Delta_{\mathcal H,Y}$ be the distributive constant of the family $\mathcal H=\{H_{i_j}\}_{j=1}^\ell$ with respect to $Y$. Then
$$e_Y({\bf c})\ge \frac{\delta_Y}{\Delta_{\mathcal H,Y}}(c_{i_1}+\cdots+c_{i_\ell}).$$
\end{lemma}

The following theorem is due to M. Ru and N. Sibony \cite{RS}.
\begin{theorem}[{reformulation of \cite[Theorem 4.8]{RS}}]\label{2.3}
Let $f $ be a linearly non-degenerate holomorphic map from $\Delta (R)\ (0<R\le +\infty)$ into $\P^N(\C)$. Let $H_1,\ldots,H_q$ be $q$ arbitrary hyperplanes in $\P^n(\C)$. Then, for every $\varepsilon >0$, we have
\begin{align*}
\biggl\|\ &\int_0^{2\pi}\max_K\log\sum_{j\in K}\frac{\|\bbf\|}{|H_j(\bbf)|}\frac{d\theta}{2\pi}+N_W(r)\\
&\le (n+1)T_f(r)+\frac{n(n+1)}{2}(c_f+\epsilon)T_f(r),
\end{align*}
where $W=\det (f_i^{(k)}; 0\le i,k\le n)$ for a reduced representation $\bbf=(f_0,\ldots,f_N)$ of $f$. 
\end{theorem}
Note that, in the original theorem \cite[Theorem 4.8]{RS}, the last term of the right hand side of the above inequality is more complex and the inequality holds for all $r$ outside an exceptional set  $S\subset(0,R)$ such that $\int_S{\rm exp}((c_f+\epsilon)T_f(r))dr<+\infty$. In this reformulation, we just simplify that term but the exceptional set is a subset $S'\subset(0,R)$ such that $\int_{S'}{\rm exp}((c_f+\epsilon')T_f(r))dr<+\infty$ for some $\epsilon'>0$.

Let $\mathcal Q=\{Q_1,\ldots,Q_q\}$ be a family moving hypersurfaces in $\P^N(\C)$ given by
$$ Q_i(z)({\bf x})=\sum_{I\in\mathcal T_{d_i}}a_{iI}(z){\bf x}^I,$$
where ${\bf x}=(x_0,\ldots,x_N)$, ${\bf x}^I=x_0^{i_0}\cdots x_N^{i_N}$ for $I=(i_0,\ldots,i_N)$. Denote by $\mathcal C_{\mathcal Q}$ the set of all non-negative functions $h : \mathbf{C}^m\setminus A\longrightarrow [0,+\infty]$, which are of the form
$$ h=\dfrac{|g_1|+\cdots +|g_l|}{|g_{l+1}|+\cdots +|g_{l+k}|}, $$
where $k,l\in\N,\ g_1,...., g_{l+k}\in\mathcal K_{\mathcal Q}\setminus\{0\}$ and $A$ is a discrete subset of $\Delta(R)$, which may depend on
$g_1,....,g_{l+k}$.  Then, for $h\in\mathcal C_{\mathcal Q}$ we have
$$\int_{0}^{2\pi}\log^+(h)\frac{d\theta}{2\pi}=O(\max T_{a_{iI}/a_{iJ}}(r)).$$
Also, for every moving hypersurface $Q$ in $\mathcal K_{\mathcal Q}[x_0,\ldots,x_N]$ of degree $d$, we have
$$ Q(z)({\bf x})\le c(z)\|{\bf x}\|^d $$
for some $c\in\mathcal C_{\mathcal Q}$.

\begin{lemma}[{see \cite[Lemma 3.2]{Q22c}}]\label{2.4}
Let $V$ be a projective variety of $\P^N(\C)$. With the above notation, let $1\le j_1\le\cdots\le j_k\le q$. Suppose that there exists $z_0\in\Delta(R)$ such that $V\cap\bigcap_{s=1}^kQ_{j_s}(z_0)^*=\emptyset.$ Then we have $V\cap\bigcap_{s=1}^kQ_{j_s}(z)^*=\emptyset$ for every $z\in\Delta(R)$ outside a discrete subset, and there exists a function $c\in\mathcal C_{\mathcal K}$ such that
$$ \|{\bf f}(z)\|\le c(z)\max_{1\le s\le k}\{Q_{j_s}({\bf f})(z)\}.$$
\end{lemma}

\section{Proof of Theorem \ref{1.1}}
Replacing $Q_{j}$ by $Q^{\frac{d}{d_{j}}}$ if necessary, we may assume that $Q_{1},\ldots,Q_{q}$ have the same degree $d$ and $Q_i=\sum_{I\in\mathcal T_d}a_{iI}x^I\ (i=1,\ldots,q)$. Take a point $z_0$  such that $a_{iI_0}(z_0)\ne 0$ for all $i$, and 
$$\Delta_{V}=\max_{\Gamma\subset\{1,\ldots,q\}}\dfrac{\sharp\Gamma}{\dim V-\dim V\cap\bigcap_{j\in\Gamma}\tilde Q_i(z_0)^*}.$$
It is suffice for us to consider the case where $q>\Delta_V(n+1)$. Denote by $\sigma_1,\ldots,\sigma_{n_0}$ all bijections from $\{0,\ldots,q-1\}$ into $\{1,\ldots,q\}$, where $n_0=q!$. For each $\sigma_i$, it is easy to see that $\bigcap_{j=0}^{q-2}\tilde Q_{\sigma_i(j)}(z_0)^*\cap V=\emptyset$. Then there exists a smallest index $\ell_i\le q-2$ such that $\bigcap_{j=0}^{\ell_i}\tilde Q_{\sigma_i(j)}(z_0)^*\cap V=\emptyset$. Hence $\bigcap_{j=0}^{\ell_i}\tilde Q_{\sigma_i(j)}(z)^*\cap V=\emptyset$ for generic points $z$ and for all $i=1,\ldots,n_0$. Denote by $\mathcal S$ the set of all points $z\in\Delta(R)$ such that $\bigcap_{j=0}^{\ell_i}\tilde Q_{\sigma_i(j)}(z)^*\cap V\ne \emptyset$ for some $i$. Then $\mathcal S$ is a discrete subset of $\Delta(R)$.

By Lemma \ref{2.4}, there is a function $A\in\mathcal C_{\mathcal Q}$, chosen common for all $\sigma_i$, such that
$$ \|\bbf (z)\|^d\le A(z)\max_{0\le j\le \ell_i}\frac{|\tilde Q_{\sigma_i(j)}(\bbf)(z)|}{\|\tilde Q_{\sigma_i(j)}(z)\|}\ \forall i=1,\ldots,n_0.$$
Denote by $S(i)$ the set of all $z$ not in $\mathcal S$ such that $\tilde Q_j(\bbf)(z)\ne 0$ for all $j=1,\ldots,q$ and
$$ \frac{|\tilde Q_{\sigma_i(0)}(\bbf)(z)|}{\|\tilde Q_{\sigma_i(0)}(z)\|}\le \frac{|\tilde Q_{\sigma_i(1)}(\bbf)(z)|}{\|\tilde Q_{\sigma_i(1)}(z)\|}\le\cdots\le \frac{|\tilde Q_{\sigma_i(q-1)}(\bbf)(z)|}{\|\tilde Q_{\sigma_i(q-1)}(z)\|}.$$
Therefore, for every generic point $z\in S(i)$, we have
$$\prod_{j=1}^q\dfrac{\|\bbf (z)\|^d\|\tilde Q_{j}(z)\|}{|\tilde Q_j(\bbf)(z)|}\le  C(z)\prod_{j=0}^{\ell_j}\dfrac{\|\bbf (z)\|^d\|\tilde Q_{\sigma_i(j)}(z)\|}{|\tilde Q_{\sigma_i(j)}(\bbf)(z)|},$$
where $C(z)=\sum_{i=1}^{n_0}A(z)^{q-\ell_i-1}\in\mathcal C_{\mathcal Q}$.

For $z\not\in\mathcal S$, consider the mapping $\Phi_z$ from $V$ into $\P^{q-1}(\C)$ defined by
$$\Phi_z({\bf x})=(\tilde Q_1(z)(x):\cdots : \tilde Q_{q}(z)(x))$$
for every ${\bf x}=(x_0:\cdots:x_N)\in V$, where $x=(x_0,\ldots,x_N)$. We set
$$\tilde\Phi_z(x)=(\tilde Q_1(z)(x),\ldots ,\tilde Q_{q}(z)(x)).$$
Let $Y_z=\Phi_z(V)$. Since $V\cap\bigcap_{j=1}^{q}\tilde Q_j(z)^*=\emptyset$, $\Phi_z$ is a finite morphism on $V$ and $Y_z$ is a projective subvariety of $\P^{q-1}(\C)$ with $\dim Y_z=n$ and of degree
$$\delta_z:=\deg Y_z\le d^{n}.\deg V=\delta.$$ 
For every ${\bf a} = (a_1,\ldots,a_q)\in\mathbb Z^q_{\ge 0}$ and ${\bf y} = (y_1,\ldots,y_q)$ we denote ${\bf y}^{\bf a} = y_{1}^{a_{1}}\ldots y_{q}^{a_{q}}$. Let $u$ be a positive integer. We set $\xi_u:=\binom{q+u}{u}$ and define the $\C$-vector space
$$ Y_{z,u}:=\C[y_1,\ldots,y_q]_u/(I_{Y_z})_u.$$
Denote by $(I_Y)_u$ the subspace of the $\mathcal K_{\mathcal Q}$-vector space $\mathcal K_{\mathcal Q}[y_1,\ldots,y_q]_u$ consisting of all homogeneous polynomials $P\in \mathcal K_{\mathcal Q}[y_1,\ldots,y_q]_u$ (including the zero polynomial) such that 
$$P(z)(\Phi_z(\bbf(z)))\equiv 0.$$
Let $(\tilde R_1,\ldots,\tilde R_p)$ be an $\mathcal K_{\mathcal Q}$-basis of $(I_Y)_u$. By enlarging $\mathcal S$ if necessary, we may assume that all zeros and poles of all nonzero coefficients of $\tilde R_i\ (1\le i\le p)$ are contained in $\mathcal S$, also all above assertions for generic points $z$ still hold for all $z\not\in\mathcal S$. Choose $\xi_u-p$ nonzero monic monomial $v_1,\ldots,v_{\xi_u-p}$ of degree of $u$ in variables $y_1,\ldots,y_q$  such that $\{\tilde R_1,\ldots,\tilde R_p,v_1,\ldots,v_{\xi_u-p}\}$ is a $\mathcal K_{\mathcal Q}$-basis of $\mathcal K_{\mathcal Q}[y_1,\ldots,y_q]_u$. 

Denote by $\mathcal T=\{T_1,\ldots,T_{\xi_u}\}$ the set of all nonzero monic monomials of degree of $u$ in variables $y_1,\ldots,y_q$. Then $\{T_1,\ldots,T_{\xi_u}\}$ is a $\mathcal K_{\mathcal Q}$-basis of $\mathcal K_{\mathcal Q}[y_1,\ldots,y_q]_u$, and also is an $\C$-basis of $\C[y_1,\ldots,y_q]_u$.

From \cite[Claim 4.3]{Q22c}, we have the following claim.
\begin{claim}\label{3.1}
There is a discrete subset $\mathcal S'$ of $\Delta (R)$ such that for all $z\not\in\mathcal S'$, we have: 
\begin{itemize}
\item[(i)] the family of equivalent classes of $v_1,\ldots,v_{\xi_u-p}$ is a basis of $Y_{z,u}$ and the family $\{\tilde R_1(z),\ldots,\tilde R_p(z)\}$ is a $\C$-basis of $(I_{Y_z})_u$;
\item[(ii)] for a subset $\{v_1',\ldots,v'_{\xi_u-p}\}$ of $\mathcal T$, if $\{\tilde R_1,\ldots,\tilde R_p,v_1',\ldots,v'_{\xi_u-p}\}$ is a $\mathcal K_{\mathcal Q}$-basis of $\mathcal K_{\mathcal Q}[y_1,\ldots,y_q]_u$ then the set of equivalent classes of $v'_1,\ldots,v'_{\xi_u-p}$ modulo $(I_{Y_z})_u$ is a $\C$-basis of $Y_{z,u}$ for everey $z\not\in\mathcal S$;
\item[(iii)] otherwise if $\{\tilde R_1,\ldots,\tilde R_p,v_1',\ldots,v'_{\xi_u-p}\}$ is linearly dependent over $\mathcal K_{\mathcal Q}$ then the set of equivalent classes of $v_1',\ldots,v'_{\xi_u-p}$ modulo $(I_{Y_z})_u$ is not a $\C$-basis of $Y_{z,u}$.
\end{itemize}
\end{claim}
\noindent
Then, we have $\xi_u-p=H_{Y_z}(u)$ for all $z$ outside $\mathcal S\cup\mathcal S'.$ Now, consider the holomorphic map $F$ from $\Delta(R)$ into $\P^{\xi_u-p-1}(\C)$ with the representation
$$ \bbF=(v_1(\tilde\Phi\circ \bbf),\ldots,v_{\xi_u-p}(\tilde\Phi\circ \bbf)). $$
Since $f$ is algebraically nondegenerate over $\mathcal K_{\mathcal Q}$, $F$ is linearly nondegenerate over $\mathcal K_{\mathcal Q}$.

Now, for $z\not\in\mathcal S\cup\mathcal S'$, we set ${\bf c}_z = (c_{1,z},\ldots,c_{q,z})\in\mathbb Z^{q},$ where
\begin{align*}
c_{i,z}:=\log\frac{\|\bbf(z)\|^d\|\tilde Q_i(z)\|}{|\tilde Q_i(\bbf)(z)|}\ge 0, \text{ for } i=1,\ldots,q.
\end{align*}
By the definition of the Hilbert weight, there are ${\bf a}_{1,z},\ldots,{\bf a}_{\xi_u-p,z}\in\mathbb N^{q}$ with
$$ {\bf a}_{i,z}=(a_{i,1,z},\ldots,a_{i,q,z}), $$
where $a_{i,j,z}\in\{1,\ldots,\xi_u\},$ such that the residue classes modulo $(I_Y)_u$ of ${\bf y}^{{\bf a}_{1,z}},\ldots,{\bf y}^{{\bf a}_{\xi_u-p,z}}$ form a basic of $\C[y_1,\ldots,y_q]_u/(I_{Y_z})_u$ and
\begin{align*}
S_Y(u,{\bf c}_z)=\sum_{i=1}^{\xi_u-p}{\bf a}_{i,z}\cdot{\bf c}_z.
\end{align*}
Note that ${\bf y}^{{\bf a}_{i,z}}\in\mathcal T$ and the set $\{\tilde R_1,\ldots,\tilde R_p,{\bf y}^{{\bf a}_{1,z}},\ldots,{\bf y}^{{\bf a}_{\xi_u-p,z}}\}$ is a basis of $\mathcal K_{\mathcal Q}[y_1,\ldots,y_q]$ (by Claim \ref{3.1}(iii)). Therefore, the set of equivalent classes of $\{{\bf y}^{{\bf a}_{1,z}},\ldots,{\bf y}^{{\bf a}_{\xi_u-p,z}}\}$ is a basis of $\dfrac{\mathcal K_{\mathcal Q}[y_1,\ldots,y_q]_u}{I(Y)_u}$. Then $ {\bf y}^{{\bf a}_{i,z}}=L_{i,z}(v_1,\ldots,v_{H_Y(u)})\ \text{modulo }I(Y)_u, $ where $L_{i,z}\ (1\le i\le \xi_u-p)$ are $\mathcal K_{\mathcal Q}$-independent linear forms with coefficients in $\mathcal K_{\mathcal Q}$.
We have
\begin{align*}
\log\prod_{i=1}^{\xi_u-p} |L_{i,z}(\bbF(z))|&=\log\prod_{i=1}^{\xi_u-p}\prod_{j=1}^q|\tilde Q_j(\bbf)(z)|^{a_{i,j,z}}\\
&= -S_Y(u,{\bf c}_z)+du(\xi_u-p)\log \|\bbf(z)\| +\log C_1(z),
\end{align*}
where $C_1\in\mathcal C_{\mathcal Q}$. Note that the number of these linear forms $L_{i,z}$ is finite, at most $\xi_u$. Denote by $\mathcal L$ the set of all $L_{i,z}$ occurring in the above inequalities. The above inequality follows that
\begin{align*}
\log\prod_{i=1}^{\xi_u-p}\dfrac{\|\bbF(z)\|\cdot \|L_{i,z}\|}{|L_{i,z}(\bbF(z))|}= &S_Y(u,{\bf c}_z)-du(\xi_u-p)\log \|\bbf(z)\| \\
&+(\xi_u-p)\log \|\bbF(z)\|+\log C_2,
\end{align*}
where $C_2(z)\in\mathcal C_{\mathcal Q}$. Then, we have
\begin{align}\label{3.2}
\begin{split}
S_Y(u,{\bf c}_z)\le&\max_{\mathcal J\subset\mathcal L}\log\prod_{L\in \mathcal J}\dfrac{\|\bbF(z)\|\cdot \|L\|}{|L(\bbf(z))|}+du(\xi_u-p)\log \|\bbf(z)\|\\
& -(\xi_u-p)\log \|\bbF(z)\|+\log C_2(z),
\end{split}
\end{align}
where the maximum is taken over all subsets $\mathcal J\subset\mathcal L$ with $\sharp\mathcal J=\xi_u-p$ and $\{L|L\in\mathcal J\}$ is linearly independent over $\mathcal K$.
From Theorem \ref{2.1} we have
\begin{align}\label{3.3}
\dfrac{1}{u(\xi_u-p)}S_{Y_z}(u,{\bf c}_z)\ge\frac{1}{(n+1)\delta_z}e_{Y_z}({\bf c}_z)-\frac{(2n+1)\delta_z}{u}\max_{1\le i\le q}c_{i,z}
\end{align}
Combining (\ref{3.2}) and (\ref{3.3}), we get
\begin{align}\label{3.4}
\begin{split}
&\frac{1}{(n+1)\delta_z}e_{Y_z}({\bf c}_z)\\
&\le\dfrac{1}{u(\xi_u-p)}\max_{\mathcal J\subset\mathcal L}\log\prod_{L\in \mathcal J}\dfrac{\|\bbF(z)\|\cdot \|L\|}{|L(\bbF(z))|}\\
&\ \ \ \ +\frac{(2n+1)\delta_z}{u}\max_{1\le i\le q}c_{i,z}+\log^+C_3(z)\\
&\le\dfrac{1}{u(\xi_u-p)}\max_{\mathcal J\subset\mathcal L}\prod_{L\in\mathcal J}\dfrac{\|\bbF(z)\|\cdot \|L\|}{|L(\bbF(z))|}\\
&\ \ \ \ +\frac{(2n+1)\delta_z}{u}\sum_{1\le i\le q}\log\frac{\|\bbf(z)\|^d\|\tilde Q_i(z)\|}{|\tilde Q_i(\bbf)(z)|}+\log^+C_3(z),
\end{split}
\end{align}
where $C_3(z)\in\mathcal C_{\mathcal Q}$, for every $z\in \Delta (R)$ outside a discrete subset.

Fix a point $z\in\Delta(R)\setminus(\mathcal S\cup \mathcal S')$. Choose $i\in\{1,\ldots,n_0\}$ such that
$$ e_{\sigma_i(0),z}\le e_{\sigma_i(1),z}\le\cdots\le e_{\sigma_i(q-1),z}.$$
Since $\bigcap_{j=0}^{\ell_i-1}\tilde Q_{\sigma_i(j)}(z)^*\cap V\ne\emptyset$, by Lemma \ref{2.2}, we have
\begin{align}\label{3.5}
\begin{split}
\Delta_{V}e_{Y_z}({\bf c}_z)&\ge (c_{\sigma_i(0),z}+\cdots +c_{\sigma_i(\ell_i),z})\cdot\delta_z\\
&=\delta_z\biggl(\sum_{j=0}^{\ell_i}\log\frac{\|\bbf(z)\|^d\|\tilde Q_{\sigma_i(j)}(z)\|}{|\tilde Q_{\sigma_i(j)}(\bbf)(z)|}\biggl ).
\end{split}
\end{align}
Then, from (\ref{3.2}), (\ref{3.4}) and (\ref{3.5}) we have
\begin{align}\label{3.6}
\begin{split}
\frac{1}{\Delta_{V}}\log &\prod_{i=1}^q\dfrac{\|\bbf (z)\|^d\|\tilde Q_i(z)\|}{|\tilde Q_i(\bbf)(z)|}\\
&\le\dfrac{n+1}{u(\xi_u-p)}\max_{\mathcal J\subset\mathcal L}\log\prod_{L\in\mathcal J}\dfrac{\|\bbF(z)\|\cdot \|L\|}{|L(\bbF(z))|}\\
&+\frac{(2n+1)(n+1)\delta_z}{u}\sum_{1\le i\le q}\log\frac{\|\bbf(z)\|^d\|\tilde Q_{i}(z)\|}{|\tilde Q_{i}(\bbf)(z)|}\\
&+\frac{1}{\Delta_{V}}\log C(z)+(n+1)\log^+C_3(z)
\end{split}
\end{align}
for every $z\in \Delta (R)$ outside a discrete subset.

Denote by $\Psi$ the set of all the coefficients of all linear forms $L_{i,z}$ and suppose that $\Psi=\{a_1,\ldots,a_{q_0}\}$. Then, we see that $\sharp\mathcal L\le\xi_u,\sharp\Psi=q_0\le \xi_u(\xi_u-p)$. For each positive integer $m$, denote by $\mathcal L(\Psi(m))$ the $C-$vector space generated by the set $\{a_1^{i_1}\ldots a_{q_0}^{i_{q_0}}|i_j\ge 0\ \text{ and }\sum_{j=1}^{q_0}i_j\le m\}$. By Remark 3.4 in \cite{TQ10}, there exists the smallest integer $p'$ such that 
$$ \frac{\dim\mathcal L(\Psi (p'+1))}{\dim\mathcal L(\Psi (p'))}\le 1+\frac{\epsilon}{2\Delta_{V}(n+1)}.$$
Put $s=\dim\mathcal L(\Psi (p')), t=\dim\mathcal L(\Psi (p'+1))$. Again, by  \cite[Remark 3.4]{TQ10}, we have
\begin{align*}
t&\le \left [\bigl(1+\frac{\epsilon}{2(n+1)\Delta_{V}}\bigl)^{\bigl[\frac{\sharp\Psi}{\log^2(1+\frac{\epsilon}{2(n+1)\Delta_{V}})}\bigl]+1}\right ]\\
&\le \left[\bigl(1+\frac{\epsilon}{2(n+1)\Delta_{V}}\bigl)^{\bigl[\frac{d^n\deg V(u+1)^{n+q}}{\log^2(1+\frac{\epsilon}{2(n+1)\Delta_{V}})}\bigl]+1}\right].
\end{align*}
Here, the last inequality comes from the fact that $\xi_u\le (u+1)^q$ and 
$$\xi_u-p\le\delta\binom{n+u}{n}\le d^n\deg V\binom{n+u}{n}\le d^n\deg V(u+1)^n.$$

Choose $\{b_1, \ldots ,b_s\}$ an $\C$-basis of $\mathcal L(\Psi(p'))$ and $\{b_1, \ldots ,b_t\}$ an $\C$-basis of $\mathcal L(\Psi(p'+1))$.
Consider the holomorphic map $\tilde F:\Delta(R)\rightarrow\mathbb{P}^{t(\xi_u-p)-1}(\mathbb{C})$ with a presentation 
$$\tilde\bbF=(b_1v_1(\tilde\Phi\circ \bbf),\ldots,b_1v_{\xi_u-p}(\tilde\Phi\circ \bbf),\ldots, b_tv_1(\tilde\Phi\circ \bbf),\ldots,b_tv_{\xi_u-p}(\tilde\Phi\circ \bbf)).$$
Note that $\|\ T_{\tilde F}(r)=duT_f(r)+o(T_f(r))$ and $c_{\tilde F}=\frac{1}{du}c_f$. By Theorem \ref{2.3}, we have
\begin{align}\label{3.7}
\begin{split}
&\biggl \|\ s\int_0^{2\pi}\max_{\mathcal J\subset\mathcal L}\log\prod_{L\in\mathcal J}\frac{\|\bbF\|}{|L(\bbF)|}\frac{d\theta}{2\pi}-N_{W(\tilde F)}(r)\\
&\le \int_0^{2\pi}\max_{\mathcal J\subset\mathcal L}\log\prod_{L\in\mathcal J}\prod_{i=1}^s\frac{\|\tilde\bbF\|}{|b_iL(\bbF)|}-N_{W(\tilde F)}(r)+o(T_f(r))\\
&\le t(\xi_u-p)udT_f(r)+\frac{(t(\xi_u-p)-1)t(\xi_u-p)}{2}(c_{\tilde F}+\frac{\epsilon'}{2du})T_{\tilde F}(r),
\end{split}
\end{align}
where $\max_{\mathcal J\subset\mathcal L}$ is taken over all subsets $\mathcal J$ of the system $\mathcal L$ of linear forms such that $\mathcal J$ is linearly independent over $\mathcal {K}_{\mathcal {Q}}$, $\epsilon'$ is an arbitrary positive number. Then, by integrating (\ref{3.6}) and using the above inequality, we obtain
\begin{align*}
\biggl\|\ &\left(\frac{1}{\Delta_{V}}-\frac{(2n+1)(n+1)\delta}{u}\right)\sum_{i=1}^qm_f(r,Q_i)\\
&\le\frac{d(n+1)t}{s}T_f(r)-\frac{(n+1)}{u(\xi_u-p)s}N_{W(\tilde F)}(r)\\
&+\frac{(t(\xi_u-p)-1)t(n+1)}{2su}\left(c_{\tilde F}+\frac{\epsilon'}{2du}\right)T_{\tilde F}(r).
\end{align*}

We now estimate the quantity $N_{W(\tilde F)}(r)$. Let $z\in\Delta(R)$ which is neither zero nor pole of any coefficients of $\tilde Q_i\ (1\le i\le q)$ and $b_i\ (1\le i\le t)$. We set 
$$c_{i}=\max\{0,\nu^0_{\tilde Q_i(\bbf)}(z)-\xi_u+p+1\}\ (1\le i\le q) \text{ and } {\bf c}=(c_{1},\ldots,c_{q})\in\mathbb Z^q_{\ge 0}.$$
Then there are 
$${\bf a}_i=(a_{i,1},\ldots,a_{i,q}),a_{i,s}\in\{1,...,u\}$$
such that ${\bf y}^{{\bf a}_1},...,{\bf y}^{{\bf a}_{H_Y(u)}}$ is a basic of $\C[y_1,\ldots,y_q]_u/(I_{Y_z})_u$ and
$$ S_{Y_z}(u,{\bf c})=\sum_{i=1}^{H_Y(u)}{\bf a}_i\cdot{\bf c}.$$
Similarly as above, we write ${\bf y}^{{\bf a}_i}=L_i(v_1,...,v_{\xi_u-p})$, where $L_1,...,L_{\xi_u-p}$ are linearly independent linear forms. We see that 
$$ \nu^0_{W(\tilde F)}(z)\ge t\sum_{i=1}^{H_Y(u)}\max\{0,\nu^0_{L_i(\bbF)}(z)-n_u\},$$
where $n_u=(\xi_u-p)t-1$. It is easy to see that
$$ \nu^0_{L_i(\bbF)}(z)=\sum_{j=1}^qa_{i,j}\nu^0_{\tilde Q_j(\bbf)}(z),$$
and hence
$$ \max\{0,\nu^0_{L_i(\bbF)}(z)-n_u\}\ge\sum_{j=1}^qa_{i,j}c_{j}={{\bf a}_i}\cdot{\bf c}. $$
Thus, we have
\begin{align}\label{3.8}
\nu^0_{W(\tilde F)}(z)\ge t\sum_{i=1}^{H_Y(u)}{{\bf a}_i}\cdot{\bf c}=tS_Y(u,{\bf c}).
\end{align}
Choose an index $\sigma_{i_0}$ such that $\nu^0_{\tilde Q_{\sigma_{i_0}(0)}(\bbf)}(z)\ge \nu^0_{\tilde Q_{\sigma_{i_0}(1)}(\bbf)}(z)\ge\cdots\ge \nu^0_{\tilde Q_{\sigma_{i_0}(q-1)}(\bbf)}(z)$. By Lemma \ref{2.2} we have
\begin{align*}
\Delta_{V}e_{Y_z}({\bf c})&\ge (c_{\sigma_{i_0}(0),z}+\cdots +c_{\sigma_{i_0}(\ell_{i_0}),z})\cdot\delta_z\\
&=\delta_z\cdot\sum_{j=1}^{l_{i_0}}\max\{0,\nu^0_{\tilde Q_{\sigma_{i_0}(j)}(\bbf)}(z)-n_u\}\\
&=\delta_z\cdot\sum_{j=1}^{q}\max\{0,\nu^0_{\tilde Q_j(\bbf)}(z)-n_u\}+O(\nu_{R^{i_0}}(z)),
\end{align*}
where $R^{i}$ is the resultant of the family $\{Q_{\sigma_i(j)}\}_{j=0}^{\ell_i}$ for $i=1,\ldots,q$. 

On the other hand, by Theorem \ref{2.1} we have that 
\begin{align*}
 S_{Y_z}(u,{\bf c}) &\ge\frac{u(\xi_u-p)}{(n+1)\delta_z}e_Y({\bf c})-(2n+1)\delta_z(\xi_u-p)\max_{1\le i\le q}c_{i}+O(\nu_{R^{i_0}}(z))\\
&\ge\left(\frac{u(\xi_u-p)}{\Delta_{V}(n+1)}-(2n+1)\delta_z (\xi_u-p)\right)\\
&\ \ \times \sum_{j=1}^{q}\max\{0,\nu^0_{\tilde Q_j(\bbf)}(z)-n_u\}+O(\nu_{R^{i_0}}(z)).
\end{align*}
Combining this inequality and (\ref{3.8}), we have
\begin{align*}
\dfrac{(n+1)}{u(\xi_u-p)s}\nu^0_{W(\tilde F)}(z)&\ge\dfrac{t}{us}\left(\frac{u}{\Delta_{V}}-(2n+1)(n+1)\delta_z\right )\\
&\times\sum_{j=1}^{q}\max\{0,\nu^0_{\tilde Q_j(\bbf)}(z)-n_u\}+O(\sum_{i=1}^{n_0}\nu_{R^{i}}(z)).
\end{align*}
Integrating both sides of this inequality, we obtain 
\begin{align*}
\biggl \|\ \dfrac{(n+1)}{u(\xi_u-p)s}N_{W(\tilde F)}(r)&\ge \dfrac{t}{s}\left(\frac{1}{\Delta_{V}}-\frac{(2n+1)(n+1)\delta}{u}\right )\\
&\ \ \times\sum_{j=1}^{q}\left(N_{Q_j(f)}(r)-N^{[n_u]}_{Q_j(f)}(r)\right)+o(T_f(r)).
\end{align*}
Seting $m_0=\frac{1}{\Delta_{V}}-\frac{(2n+1)(n+1)\delta}{u}$ and combining inequalities (\ref{3.7}) with the above inequality, we get
\begin{align*}
\biggl\|\ \sum_{i=1}^{q}m_f(r,Q_i)\le &\frac{d(n+1)t}{sm_0}T_f(r)-\frac{t}{s}\sum_{j=1}^{q}\left(N_{Q_i(f)}(r)-N^{[n_u]}_{Q_j(f)}(r)\right)\\
&+\frac{(t(\xi_u-p)-1)t(n+1)}{2sm_0u}(c_{\tilde F}+\frac{\epsilon'}{2du})T_{\tilde F}(r)+o(T_f(r))
\end{align*}
This inequality implies that
\begin{align}\label{3.9}
\begin{split}
\bigl\|\ &\bigl(q-\frac{(n+1)t}{sm_0}\bigl)T_f(r)\\
&\le\sum_{j=1}^{q}\frac{1}{d}N^{[n_u]}_{Q_j(f)}(r)+\frac{(t(\xi_u-p)-1)t(n+1)}{2sdm_0u}(c_{f}+\epsilon')T_{f}(r).
\end{split}
\end{align}

\noindent
a) We choose $u=\lceil 2\Delta_{V}(2n+1)(n+1)d^n\deg V(\Delta_{V}(n+1)+\epsilon)\epsilon^{-1}\rceil$. 

Then $u\ge \biggl\lceil \dfrac{\Delta_{V}(2n+1)(n+1)\delta(\Delta_{V}(n+1)+\epsilon)}{\Delta_{V}(n+1)+\epsilon-\frac{\Delta_{V}(n+1)t}{s}}\biggl\rceil,$ and we have:
\begin{align*}
\bullet\ &q-\frac{(n+1)t}{sm_0}\ge q-\frac{\Delta_V(n+1)t}{(1-\Delta_V(2n+1)(n+1)\delta/u)s}\\
&\ge q-\Delta_{V}(n+1)-\epsilon;\\
\bullet\ &n_u+1=(\xi_u-p)t\\
&\le  d^n\deg V(u+1)^n\left[\bigl(1+\frac{\epsilon}{2(n+1)\Delta_{V}}\bigl)^{\bigl[\frac{d^n\deg V(u+1)^{n+q}}{\log^2(1+\frac{\epsilon}{2(n+1)\Delta_{V}})}\bigl]+1}\right]=L;\\
\bullet\ &\frac{(t(\xi_u-p)-1)t(n+1)}{2sdm_0u}< \frac{t(n+1)}{sm_0}\cdot\frac{(L-1)}{2du}\\
&\le \frac{(\Delta_{V}(n+1)+\epsilon)(L-1)}{2du}.
\end{align*}
Then, from (\ref{3.9}) we have
\begin{align*}
\bigl\|\ &(q-\Delta_{V}(n+1)-\epsilon)T_f(r)\\
&\le\sum_{j=1}^{q}\frac{1}{d}N^{[L-1]}_{Q_j(f)}(r)+\frac{(\Delta_{V}(n+1)+\epsilon)(c_f+\epsilon')(L-1)}{2du}T_f (r).
\end{align*}
The assertion a) is proved.
 
\noindent
b) If all $Q_i$ are fixed hypersurfaces then $t=s=1$. Choosing $u'=\lceil \Delta_{V}(2n+1)(n+1)d^n\deg V(\Delta_{V}(n+1)+\epsilon)\epsilon^{-1}\rceil$ and replacing $u$ in the above by $u'$, we have
\begin{align*}
\bullet\ &q-\frac{(n+1)}{m_0}\ge q-\Delta_{V}(n+1)-\epsilon,\\
\bullet\ &n_{u'}+1=\xi_{u'}-p<\delta\binom{n+u'}{n}\le d^n\deg V\binom{n+u'}{n}\\
&\le \left[d^n\deg Ve^n\left(1+\frac{u'}{n}\right)^{n}\right]\\
&\le \left[d^{n^2+n}(\deg V)^{n+1}e^n(2n+5)^n(\Delta^2_V(n+1)\epsilon^{-1}+\Delta_V)^n\right]=L'.
\end{align*}
Similarly as above, from (\ref{3.9}) we have the desired inequality of the assertion b).

Hence, the proof of the theorem is completed.\hfill$\square$

\vskip0.1cm
\noindent
\textbf{\textit{Remark:}} For the case $\Delta(R)=\C$, we have $c_f=0$ if $f$ is nonconstant. In this case, (\ref{3.7}) can be re-written as follows:
\begin{align}\label{3.10}
\begin{split}
\biggl \|\ s&\int_0^{2\pi}\max_{\mathcal J\subset\mathcal L}\log\prod_{L\in\mathcal J}\frac{\|\bbF\|}{|L(\bbF)|}\frac{d\theta}{2\pi}-N_{W(\tilde F)}(r)\\
&\le \int_0^{2\pi}\max_{\mathcal J\subset\mathcal L}\log\prod_{L\in\mathcal J}\prod_{i=1}^s\frac{\|\tilde\bbF\|}{|b_iL(\bbF)|}-N_{W(\tilde F)}(r)+o(T_f(r))\\
&\le t(\xi_u-p)udT_f(r)+o(T_f(r)).
\end{split}
\end{align}
Here, the notation $``\|''$ means that the inequality hold for all $r\in [1,+\infty)$ outside a set of finite measure.

Using (\ref{3.10}) instead of (\ref{3.7}), from the above proof, we obtain the following second main theorem for holomorphic maps from $\C$.
\begin{theorem}\label{3.11}
 Let $V$ be a smooth subvriety of dimension $n$ of $\mathbb{P}^{N}(\mathbb{C})$. Let $f:\C\rightarrow V$ be a holomorphic mapping. Let $\mathcal {Q}=\{Q_{1},\ldots,Q_{q}\}$ be a set of slowly (with respect to $f$) moving hypersurfaces with the distributive constant $\Delta_{V}$ with respect to $V$. Assume that $f$ is algebraically non-degenerate over $\mathcal{K}_{\mathcal {Q}}$. Let $d=lcm(\deg Q_1,\ldots,\deg Q_q)$. 
Then for every $(n+1)\Delta_V>\epsilon>0,$
\begin{align}\label{3.12}
\bigl\|\ (q-\Delta_{V}(n+1)-\epsilon)T_f(r)\le\sum_{j=1}^{q}\frac{1}{\deg Q_j}N^{[L]}_{f}(r,Q_j), 
\end{align}
where 
$$L=d^n\deg V(u+1)^n\left[\bigl(1+\frac{\epsilon}{2(n+1)\Delta_{V}}\bigl)^{\bigl[\frac{d^n\deg V(u+1)^{n+q}}{\log^2(1+\frac{\epsilon}{2(n+1)\Delta_{V}})}\bigl]+1}\right],$$
with $u=\lceil 2\Delta_{V}(2n+1)(n+1)d^n\deg V(\Delta_{V}(n+1)+\epsilon)\epsilon^{-1}\rceil$.

Moreover, if all $Q_i\ (1\le i\le q)$ are assumed to be fixed hypersurfaces, then for every $\epsilon >0$ we have the inequality $(\ref{3.12})$ with
$L=\bigl[d^{n^2+n}(\deg V)^{n+1}e^n(2n+5)^n(\Delta^2_V(n+1)\epsilon^{-1}+\Delta_V)^n\bigl]$.
\end{theorem}
We also note that our proof is valid for the case of holomorphic maps from higher dimension complex spaces $\C^m$ into $V$. This theorem is an improvement of many previous second main theorem for hypersurface targets, such as \cite{DT20,LG,LS,Q22a,Q22c} and \cite{Ru09}.

\section*{Disclosure statement}
No potential conflict of interest was reported by the author(s).

\end{document}